\documentclass[review]{elsarticle}
\usepackage{amsfonts,amssymb,amsmath,amsbsy,amsthm,amscd}
\usepackage{lineno,hyperref}
\modulolinenumbers[5]
\usepackage[english]{babel}
\usepackage{graphicx}

\newtheorem{theorem}{Theorem}[section]

\newtheorem{remark}{Remark}

\begin{document}

\begin{frontmatter}

\title{Suppression of singularities of solutions of the Euler-Poisson system with  density-dependent damping
}


\author{Olga S. Rozanova}
\ead{rozanova@mech.math.msu.su} 

\address[1]{Department of Mechanics and Mathematrics, Moscow
State University, Moscow 119991 Russia}


\begin{abstract}
We find a sharp condition on the density-dependent coefficient of
damping of a one-dimensional repulsive Euler-Poisson system, which
makes it possible to suppress the formation of singularities in the
solution of the Cauchy problem with arbitrary smooth data. In the
context of plasma physics, this means the possibility of suppressing
the breakdown of arbitrary oscillations of cold plasma.
\end{abstract}


\def\sign{\mathop{\rm sgn}\nolimits}






\end{frontmatter}

\section{Introduction}

We consider the system of Euler-Poisson equations describing the behavior  of cold plasma 
for velocity $V$, electron density $n> 0$ and electric field
potential
  $\Psi$ in the following form:
\begin{equation}\label{1}
\dfrac{\partial V}{\partial t}+ (V \cdot \nabla ) V =-\nabla \Psi
-\nu V,\quad \dfrac{\partial n}{\partial t}+ {\rm div}\,( n V)
=0,\quad \Delta \Psi =1 -n.
\end{equation}
All components of the solution are assumed to be functions of time
$t\ge 0$ and point $x\in {\mathbb R}^n$, $\int\limits_{{\mathbb
R}^n} (n-1) \,dx = \rm const$, $\nu\ge 0$ is the  damping factor.

This system is a pressureless variant of the general Euler-Poisson
system, having numerous physical applications, see \cite{ELT} for
references. One of the crucial questions is the study of the Cauchy
problem and the analysis of the possibility of the existence of a
globally in time smooth solution. The model without pressure is
somewhat simpler from a mathematical point of view, since it allows
one to obtain criteria for the formation of a singularity from the
initial data.

In \cite {ELT}, many versions of the model without pressure,
including those with constant damping and viscosity, with both zero
and non-zero backgrounds, have been studied. In all these cases, it
is possible to find the initial data leading to a blow-up in a
finite time. Moreover, this possibility still remains if pressure
and heat diffusion are added to the  model   \cite {ChenWang}.

In recent years, the pressureless model has attracted great
interest, since it is very convenient to describe the wake wave in
the cold plasma generated by a laser pulse in order to create a new
type of accelerator, \cite{Ch_book} and references therein. It is
generally known that the plasma oscillation tends to blow-up,
forming a gradient catastrophe in the velocity component and a delta
singularity in the density component. After the moment of the
singularity formation, the cold plasma model loses its relevance;
therefore, the conditions on the initial data or other parameters
that make it possible to maintain a smooth solution as long as
possible or, possibly, guaranteeing a global in time smooth solution
is a key question for all the theory.

In plasma physics, the damping factor $\nu$ corresponds to the
frequency of the electron-ion collisions, this value is very small
from the physical point of view. Depending on the model,
 the
electron-ion collisions either can be neglected or taken into
account. In the recent paper \cite{brodin}, for a very particular
solution in the $1D$ case the authors showed numerically that if
$\nu=\nu_0 n $, where $\nu_0$  is a positive constant, than the
oscillations never blow up.

Our main question is whether this result is valid for all possible
initial data and can it be substantiated analytically? In addition,
if $ \nu $ is a smooth function of $ n $, what conditions must be
imposed on $ \nu (n) $ to ensure the global in time smoothness of
the solution to the Cauchy problem  for any given data?

In this paper, we focus on the $ 1D $ case, since explicit
analytical results can be obtained here, and we rewrite \eqref {1}
in a form more accepted in plasma physics:
\begin{equation}\label{2}
\dfrac{\partial V}{\partial t}+ V  \dfrac{\partial V}{\partial x}=-E
-\nu(n) V, \quad \dfrac{\partial E}{\partial t}+ V  \dfrac{\partial
E}{\partial x}=V, \quad n=1- \dfrac{\partial E}{\partial x},
\end{equation}
see \cite{Ch_book} for details. Here $E=\nabla \Psi$ is the vector
of electric field. System \eqref{2} will be considered together with
the Cauchy data
\begin{equation}\label{CD}
(V,E)|_{t=0}=(V_0(x), E_0(x))\in {\mathcal A} ({\mathbb R}).
\end{equation}
For $\nu\ne \rm const$ system \eqref{2}, do not belong to the
symmetric hyperbolic type, therefore we cannot guarantee that the
solution to the Cauchy problem a local solution as smooth as initial
data in the Sobolev norm. Therefore we have to prescribe the
analyticity to initial data to use the Cauchy-Kovalevskaya theorem.
to show that problem \eqref{2},
 \eqref{CD}  a local in time
unique analytical solution.

Problem \eqref{2}, \eqref{CD} was completely analyzed  for $\nu=0$
in \cite{RCh_ZAMP} and $\nu=\rm const>0$ in \cite{DRCh_AIP}, where
sharp conditions on initial data to guaranty a globally in time
smooth solution were found (see the analogous result for another
context in \cite{ELT}). It was found that even for an arbitrarily
large constant frequency of collisions there exist data implying a
finite time singularity formation.

In the present work, we show that by choosing an appropriate
density-dependent damping factor one can obtain a globally smooth
solution for any smooth Cauchy data, i.e. completely remove the
singularity formation.  In particular, for our prototypic function
$\nu(n)=\nu_0 n^\gamma$ the threshold value is $\gamma=1$. For
$\gamma> 1$ the solution to \eqref{2}, \eqref{CD} does not form the
gradient catastrophe for any choice of initial data.

The paper is organized as follows. In Section \ref{S1}, we consider
a special solution, linear with respect to the spatial variable (the
so-called affine solutions) and prove the exact condition for
eliminating blow-up. In Section \ref{S2} we prove a similar result
for arbitrary initial data. Section \ref{S4} discusses issues
related to this problem.


\section{Affine solutions}\label{S1}

First, we consider a special form of solutions:
\begin{equation}\label{as}
V=a(t) x +A(t), \quad E= b(t) x +B(t).
\end{equation}

\begin{theorem}\label{as_T} Let $f(n)\in {\mathcal A}({\mathbb R}_+) $ be a nonnegative function
satisfying condition
\begin{equation}\label{as_cond}
\int\limits_{\eta_0>0}^{+\infty}\dfrac{ f(\eta)}{\eta^2}\,
d\eta=\infty.
\end{equation}
It $\nu(n) = \epsilon f (n)$, $\epsilon={\rm const}>0$, then
derivatives of the solution  to problem \eqref{2}, \eqref{CD}, with
 data  \eqref{as} are bounded in time for any choice of the data.
Otherwise, one can find the data such that the derivatives of
solution blow up in a finite time.
\end{theorem}

\proof 
We substitute the ansatz \eqref{as} in \eqref{1} to obtain
\begin{eqnarray}\label{as_s1}
&\dot{a}=-a^2-b-\epsilon f(1-b) a,& \quad \dot{b}= (1-b) a,\\
&\dot{A}=-A (a - \epsilon f(1-b))- B, & \quad \dot{B}=(1-b) A,
\label{as_s2}
\end{eqnarray}
Since $n=1-b>0$, we consider the domain $b<1$.

 The couple of equations \eqref{as_s1} splits off from the
system, and the second couple \eqref{as_s2} is linear with respect
to $A, B$ with the coefficients found in the previous step.
Therefore, if we want to study conditions for a blowup of the
solution \eqref{as}, it is sufficient to consider the behavior of
the phase curve of the autonomous system \eqref{as_s1}, given as
\begin{eqnarray}\label{as_phase}
\dfrac {d a}{db}=   - \frac{a^2+b}{ (1-b)a}  -\epsilon
\frac{f(1-b)}{1-b}.
\end{eqnarray}
For an arbitrary $f$ equation \eqref{as_phase} cannot be integrated
explicitly, however, it can be considered as a regular perturbation
of \eqref{as_phase} at $\epsilon=0$,
 which solution is
\begin{eqnarray}\label{as_a0}
a=\pm \sqrt{1-2b+C (1-b)^2}.
\end{eqnarray}
Here the constant $C=\dfrac{a_0^2+2b_0-1}{(1-b_0)^2}$, with
$a_0=a(0)$, $b_0=b(0)$.

The analysis of the phase plane shows that  a point on the phase
plane moves from the upper half-plane $a>0$ to the lower half-plane
$a<0$ and there can come back to the upper half-plane or go to
minus-infinity. The latter signifies the blowup of derivatives of
the solution. If $C<0$, the curve on the phase plane is bounded (it
is ellipse), otherwise $a(t)$ and $b(t)$ move along a parabola
($C=0$) or hyperbola ($C>0$) and therefore go to minus-infinity
within a finite time (see \cite{RCh_ZAMP} for details).

Our main question is whether correctors due to parameter $\epsilon$
can change the behavior of trajectory going to infinity and  turn it
to the upper half-plane $a>0$.

 When analyzing the phase portrait of the perturbed system,  we point out the following elementary facts, illustrated in
 Figure 1.

 1. If the initial point of a phase curve is situated in the upper half-plane $a>0$, within a finite
 time point $(b,a)$ turns in the lower half-plane $a<0$, therefore a possible
 blowup can happen only for $a<0$  (see Figure 1, left);

 2. Since $\dfrac{d a}{d b}=  - \dfrac{a^2+b}{ (1-b)a}  -\epsilon
\dfrac{f(1-b)}{1-b}\le- \dfrac{a^2+b}{ (1-b)a}$, and $b$ decreases
with $t$ as $a<0$, then the Chaplygin theorem implies that the phase
curve of the perturbed equation, $a_\epsilon(b)$ lies higher that
the phase curve of non-perturbed equation, $a_0(b)=- \sqrt{1-2b+C
(1-b)^2}.$ Therefore for $C<0 $ the curve $a_\epsilon(b)$ always
comes back to the upper half-plane $a>0$. So, for a possible blow-up
we have to consider only the initial data corresponding to $C\ge 0$
 (see Figure 1, right).

3. Analogously, the Chaplygin theorem implies that if the data are
such that the phase curve $a_{\epsilon_1}(b)$ does not go to
infinity, then $a_{\epsilon_2}(b)$, $\epsilon_2>\epsilon_1$ does not
go to infinity as well. Therefore, we can consider for the proof
arbitrarily small $\epsilon$.

Since the perturbation by means of parameter $\epsilon$ is regular,
in a neighborhood $U_\epsilon(0)$  we can expand the solution in a
series
\begin{equation*}\label{expand}
a_\epsilon(b)= a_0(b)+\sum\limits_{k=1}^{\infty}\epsilon^k \alpha_k
(b)=a_0(b)+\epsilon \alpha_1 (b)+ o(\epsilon),
\end{equation*}
converging at any fixed $b$.

Thus, if the first corrector $\alpha_1 (b)$ is such that
$a_0(b)+\epsilon \alpha_1 (b)>0$ for some $b_*$ and arbitrary small
positive $\epsilon$, then we can guarantee that $a_\epsilon(b_*)>0$,
in other words, the trajectory came back to the upper half-plane.
The linear equations for the correctors are the following:
\begin{eqnarray}
\dfrac{d\alpha_1}{db} &=& -\alpha_1 Q -\dfrac{f(1-b)}{1-b}, \label{al1}\\
 \dfrac{d\alpha_k}{d b} &=& -\alpha_k Q +\dfrac{b \phi_k(a_0,\alpha_1,\dots,\alpha_{k-1})}{(1-b)a_0^k},\quad
k=2,\dots,\label{al_k}
\end{eqnarray}
where $Q= \dfrac{a_0^2-b}{(1-b) a_0^2}=\dfrac{1-3b+C
(1-b)^2}{(1-3b+C (1-b)^2)(1-b)} $, and $ \phi_k $ is a homogeneous
polynomial  of order $k$ from its arguments, $a_0$ is found in
\eqref{as_a0}. It can be readily found from \eqref{al1} that
\begin{equation*}\label{al1_sol}
\alpha_1(b)= \frac{(1-b)^2}{a_0(b)}
\int\limits_{b}^{b_0}\dfrac{a_0(\beta) f(1-\beta)}{(1-\beta)^3}\,
d\beta >0,\quad b<b_0=\rm const.
\end{equation*}
As for $C>0$
\begin{equation*}\label{as_asymp}
a_0(b)+\epsilon \alpha_1(b)\sim -\sqrt{C}(1-b)+ \epsilon
\frac{(1-b)^2}{a_0(b)} \int\limits_{b}^{b_0}\dfrac{
f(1-\beta)}{(1-\beta)^2}\, d\beta, \quad b\to -\infty,
\end{equation*}
then condition \eqref{as_cond} guarantees  boundedness of the phase
trajectory $a_\epsilon(b)$. If the integral \eqref{as_cond}
converges, then for sufficiently small $\epsilon$ the prevailing
term is $a_0$, and the trajectory goes to infinity. Moreover, for
sufficiently small $\epsilon$, satisfying condition
\begin{equation}\label{as_cond_gc}
-\sqrt{C}+\epsilon \int\limits_{1-b_0>0}^{+\infty}\dfrac{
f(\eta)}{\eta^2}\, d\eta<0,
\end{equation}
the solution blowup for the same initial data as in the
non-perturbed case.

For $C=0$ the analysis is analogous, but the result is different.
Namely, $a_0(b)\sim -(1-b)^{\frac{3}{2}}$ as $b\to -\infty$ and the
respective condition for  boundedness of trajectory is
\begin{equation}\label{as_cond0}
\lim\limits_{\eta\to\infty}\eta \int\limits_{\eta_0>0}^{\eta}\dfrac{
f(\tilde{\eta})}{\tilde{\eta}^{{5/2}}}\, d\tilde{\eta}=\infty.
\end{equation}

Condition \eqref{as_cond0} is predictively more mild than
\eqref{as_cond}. For example, for $f(\eta)=\eta^\gamma$,
\eqref{as_cond0} gives $\gamma>\frac12$, whereas \eqref{as_cond}
gives $\gamma\ge 1$. Nevertheless, we have to take into account the
worst situation, i.e. $C>0$.

Figure 2 shows the effect of $ \epsilon $ on the solution. It can be
seen that even for sufficiently large $ \epsilon $ the solution
first closely mimics the unperturbed case and only after some time
sharply changes its behavior (Figure 2, left). In fact, for the
threshold value $ \gamma = 1 $, the difference between the perturbed
and unperturbed cases is at first so small that it cannot be
detected numerically. Figure 2, right, shows rapidly decaying
oscillations for a sufficiently long time.

Thus, since the analytical solution to \eqref{2}, \eqref{CD} is
unique, if the data belong to the class \eqref{as}, so does the
solution. The theorem is proved.  $\Box$

\medskip
\begin{center}
\begin{figure}[h!]
\begin{minipage}{0.495\columnwidth}
\centerline{
\includegraphics[width=0.9\columnwidth]{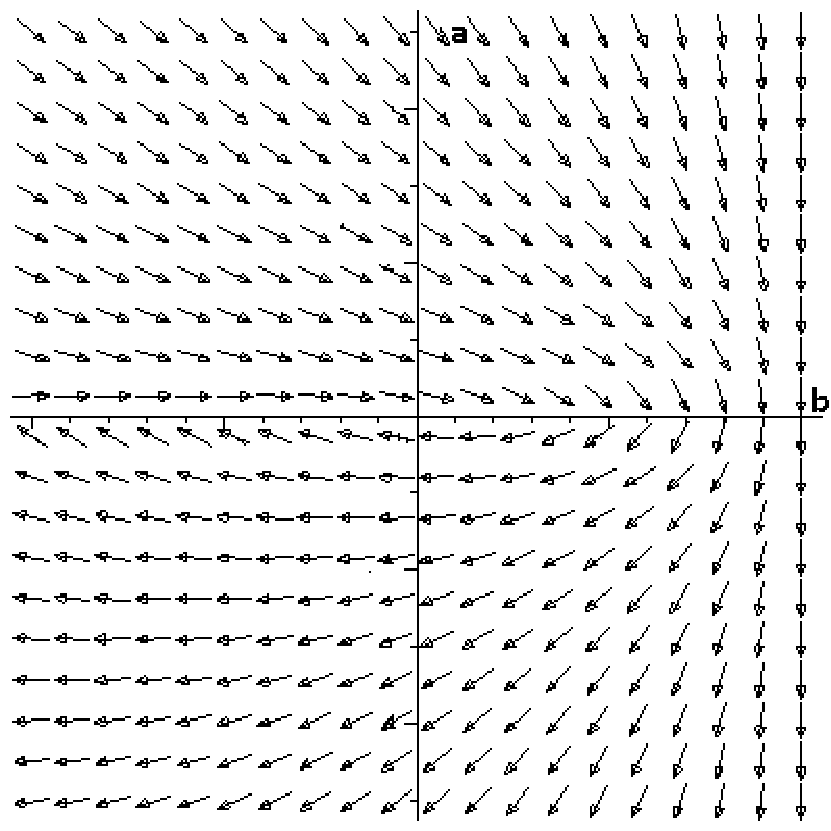}}
\end{minipage}
\begin{minipage}{0.495\columnwidth}
\centerline{
\includegraphics[width=0.9\columnwidth]{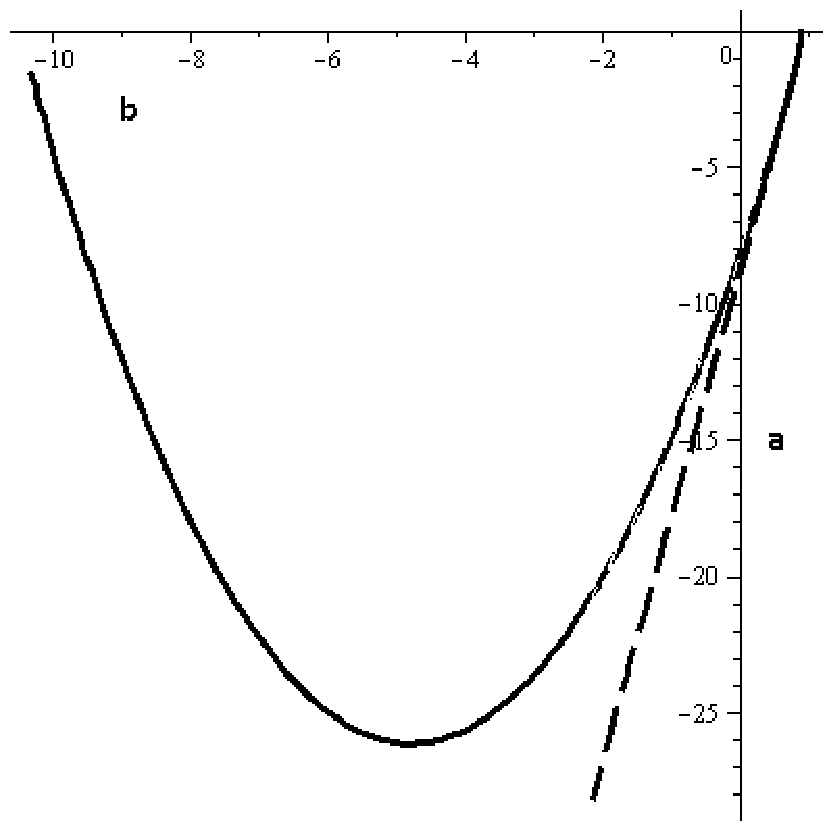}}
\end{minipage}
\caption{ $f(n)=n^2$. Left: the direction field to system
\eqref{as_s1},  $\epsilon =0.8$. Right: the phase curves starting
from the same point for \eqref{as_phase} at $\epsilon=0$,
singularity formation (dash) and $\epsilon =0.8$, smooth solution
(solid).}
\end{figure}
\end{center}

\medskip
\begin{center}
\begin{figure}[h!]
\begin{minipage}{0.495\columnwidth}
\centerline{
\includegraphics[width=0.9\columnwidth]{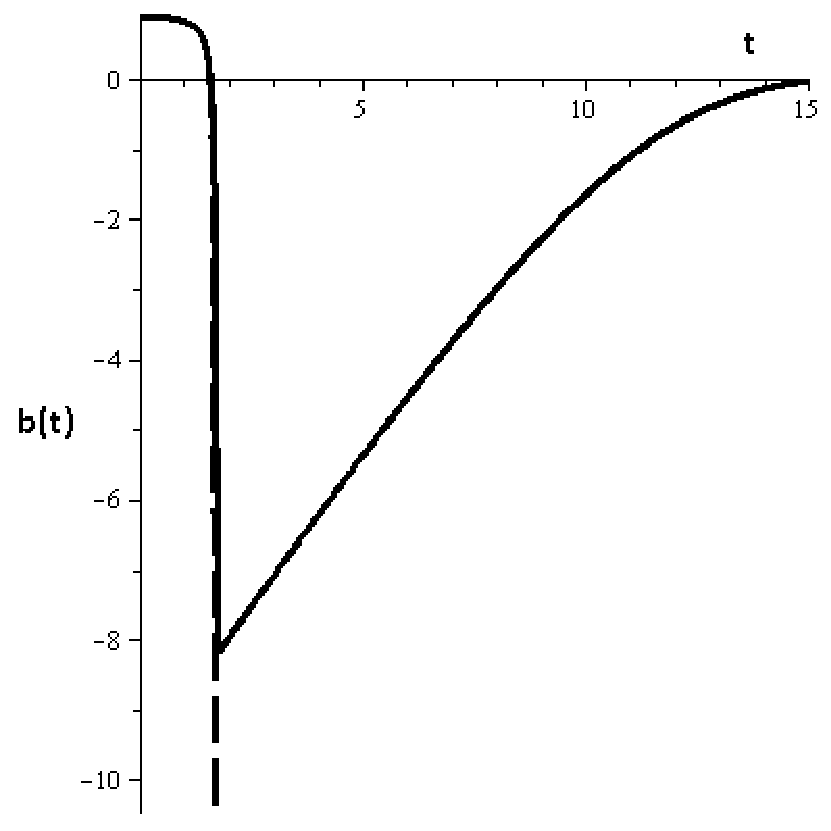}}
\end{minipage}
\begin{minipage}{0.495\columnwidth}
\centerline{
\includegraphics[width=0.9\columnwidth]{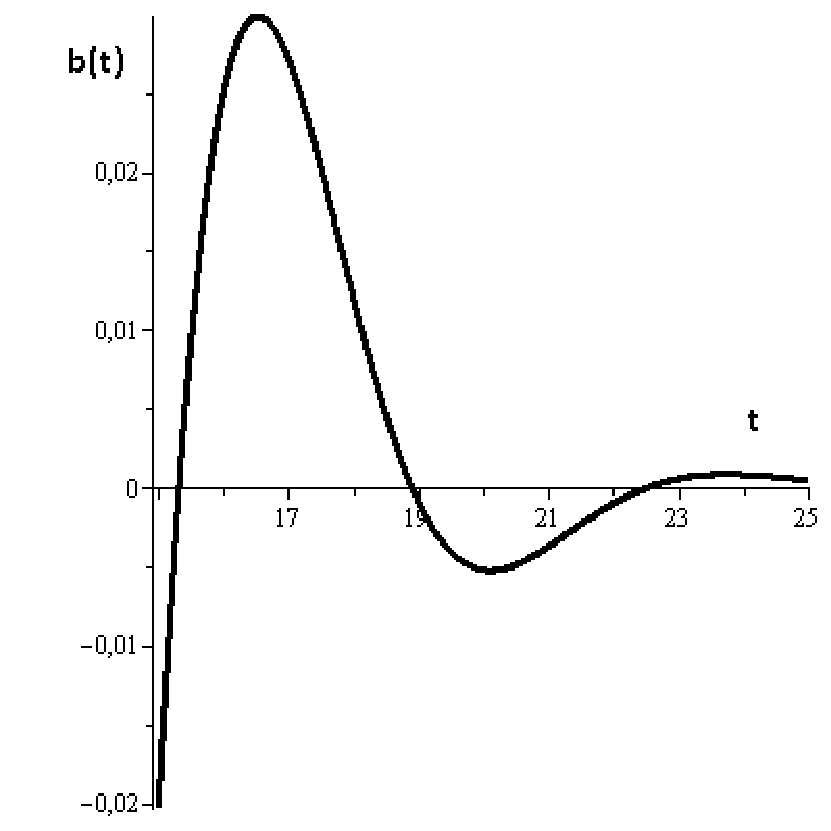}}
\end{minipage}
\caption{Solution of \eqref{as_s1} for $f(n)=n^2$. Left: the
behavior of $b(t)$
 at $\epsilon=0$, singularity formation (dash) and  $\epsilon =1$,
smooth solution (solid). Right: the behavior of $b(t)$ for $\epsilon
=1$ near equilibrium $b=0$.}
\end{figure}
\end{center}

\section{Arbitrary initial data}\label{S2}

\begin{theorem}\label{aid_T}{\bf (Main theorem)} Let
$f(n)\in {\mathcal A}({\mathbb R}_+) $ be a nonnegative function
satisfying conditions
\begin{equation}\label{aid_cond2}
    \lim\limits_{\eta\to\infty}\dfrac{\eta f'(\eta)}{f(\eta)}={\rm
    const}> 1
\end{equation}
and \eqref{as_cond}. If $\nu(n) = \epsilon f (n)$, $\epsilon =\rm
const >0$, then problem \eqref{2}, \eqref{CD} admits a global in
time classical ($C^1$-smooth) solution for any choice of the data.
Otherwise, one can find the data such that the derivatives of
solution blow up in a finite time.
\end{theorem}

\proof We denote $q=V_x$, $s=E_x$, $\xi=V_{xx}$, $\sigma=E_{xx}$ and
differentiate \eqref{2} with respect to $x$. Since $n=1-s>0$, then
it makes sense to consider only the half-plane $s<1$. Along every
characteristic line $x(t)$, starting from point $x_0\in\mathbb R$ we
get system
\begin{eqnarray}\label{aid_q}
\dot q&=& -q^2 -s - \epsilon(f(1-s) q  + V f'(1-s)\sigma  ), \\
 \dot s&=& (1-s)q,\label{aid_s}
\end{eqnarray}
complemented by initial conditions $q(0)=V_x(x_0)$, $s(0)=E_x(x_0)$.
Due to the term $V f'(1-s)\sigma$ this system is not closed. The
dynamics of $V$ can be found from  \eqref{2}:
\begin{eqnarray}
  \dot V &=& -E-\epsilon f(1-s) V, \label{aid_V} \\
  \dot E &=& V, \label{aid_E}
\end{eqnarray}
it implies
\begin{eqnarray}\label{aid_EV}
\frac{d}{dt} (V^2 +E^2)= - 2 \epsilon f(1-s) V^2\le 0,\nonumber
\end{eqnarray}
 therefore $V$ and $E$ remain bounded.

However, the equation for $\sigma$ contains $\sigma_x$ and the whole
system cannot be closed. It is the principal difficulty comparing
with the case $\epsilon=0$, treated in \cite{RCh_ZAMP}.

For the solutions \eqref{as} this problem does not arise, since
$\sigma=0$ for them.

Further, \eqref{aid_q}, \eqref{aid_s} imply
\begin{eqnarray}\label{aid_qs}
\dfrac{d q}{d s} =-\dfrac{q^2+s}{(1-s)q}-\epsilon \frac{f(1-s)}{1-s}
-\epsilon \frac{\sigma f'(1-s)}{(1-s) q} V,
\end{eqnarray}
which coincides with \eqref{as_phase}, except for the last term. We
are going to show that this term is subjected to the previous one as
$s\to -\infty$ and therefore similar to the arguments of Theorem
\ref{as_T} along every characteristic $x=x(t)$ the derivatives of
the solution are bounded. Variables $(s,q)$ correspond to $(a,b)$ in
the proof of Theorem \ref{as_T}.

Namely, to find the condition for the boundedness of $q, s$ we
consider expansion with respect to the small parameter $\epsilon$.

First of all, we introduce a new independent variable as  $s= s(t)$.
This is possible if $\dot s\ne 0$, i.e $q\ne 0$.  The blow-up
implies that $s$ tends to $-\infty$ as $t\to t_*<\infty$, for $q<0$.

Let us set
$ q(s)=q_0(s)+\epsilon q_1(s)+o(\epsilon),$ $
\sigma(s)=\sigma_0(s)+\epsilon \sigma_1(s)+o(\epsilon),$ $
\xi(s)=\xi_0(s)+\epsilon \xi_1(s)+o(\epsilon).$
Then as in \eqref{as_a0} we find
\begin{eqnarray}\label{q0(s)}
q_0(s)=\pm \sqrt{1-2s+C (1-s)^2}, \quad
C=\frac{q_0^2(0)+2s(0)-1}{(1-s(0))^2}
\end{eqnarray}
and
\begin{eqnarray}
\dfrac{dq_1}{ds} &=& -q_1 Q(s) -\dfrac{f(1-s)}{1-s}-\dfrac{\sigma_0
f'(1-s) V(s,q_0)}{(1-s)q_0}, \label{q1}
\end{eqnarray}
where $Q(s)= \dfrac{q_0^2-s}{(1-s) q_0^2}=\dfrac{1-3s+C
(1-s)^2}{(1-3s+C (1-s)^2)(1-s)} $.

To find $\sigma_0$, we get the system of linear equations
\begin{eqnarray}
  \dfrac{d \sigma_0}{ds} &=& \dfrac{(1-s)\xi_0 -2\sigma_0 q_0}{ (1-s)q_0},\label{aid_sigma0} \\
  \dfrac{d \xi_0}{d s} &=& -\dfrac{3 q_0 \xi_0+\sigma_0}{(1-s)q_0 }, \label{aid_xi0}
\end{eqnarray}

Further, taking into account \eqref{q0(s)}, from the system of
linear equations \eqref{aid_sigma0}, \eqref{aid_xi0} we have
\begin{eqnarray}\label{aid_sig}
  \sigma_0(s) &=& (s-1)^2(C_1 s +C_2 q_0(s)),
\end{eqnarray}
with constants $C_1, C_2$, depending on $s(0), q_0(0), \sigma_0(0),
\xi_0(0)$. Due to \eqref{aid_sig} and condition \eqref{aid_cond2}
the ratio $\dfrac{\sigma_0 f'(1-s)}{(1-s) q_0}$ has  the same
behavior as ${f(1-s)}$ as $s\to -\infty$.

Let us study the behavior of the term $V(s,q_0(s))$ as $s\to
-\infty$.

From \eqref{aid_q} -- \eqref{aid_E} we have
\begin{eqnarray}
  \frac{dV}{ds} &=& \frac{-E(s)-\epsilon f(1-s)V(s)}{(1-s) q(s)},\nonumber \\
  \frac{dE}{ds} &=& \frac{V(s)}{(1-s) q(s)}.\nonumber
\end{eqnarray}
Since  $V(s, q(s))=V(s, q_0(s))+O(\epsilon) $, $E(s, q(s))=E(s,
q_0(s))+O(\epsilon)$, $\epsilon\to 0$, fixed $s$, then for the zero
terms $V_0=V(s, q_0(s))$, $E_0=E(s, q_0(s))$ we obtain
\begin{eqnarray}
  \frac{dV_0}{ds} &=& \frac{-E_0(s)-\epsilon f(1-s)V_0(s)}{(1-s) q_0(s)},\label{dVs} \\
  \frac{dE_0}{ds} &=& \frac{V_0(s)}{(1-s) q_0(s)}.\label{dEs}
\end{eqnarray}

Further, for convenience, we change the independent variable once
again as
\begin{eqnarray}\label{aid_s1}
s_1(s)= - \arctan \frac{s}{\sqrt{1- 2s +C (1-s)^2}}+ \arctan
\frac{1}{\sqrt{C}}\sim \frac{1}{\sqrt{C}(1-s)},\quad  s\to
-\infty.\nonumber
\end{eqnarray}
 Then $s_1\to 0+$  and  $f(1-s) \sim f((\sqrt{C} s_1)^{-1})$ as $s\to
 -\infty$. With the new independent variable  \eqref{dVs},
 \eqref{dEs} take the form
\begin{eqnarray}
  \frac{dV_0}{ds_1} &=& E_0(s_1)+\epsilon f(1-s)V_0(s_1),\label{dVsn} \\
  \frac{dE_0}{ds_1} &=& -V_0(s_1).\label{dEsn}
\end{eqnarray}

 Let us study the structure of the solution near the point $s_1=0$,
 the hypothetic point of singularity formation.

 Further, from \eqref{dVsn} and a linear homogeneous equation
\eqref{dEsn} we have
\begin{eqnarray}\label{Vs1}
\frac{d^2 V_0 }{ds_1^2} -  \epsilon f(1-s) \frac{d V_0 }{ds_1}
+(1-\epsilon f'(1-s))V_0=0, \quad s=s(s_1),
\end{eqnarray}
where $s_1=0$ is an irregular singular point. To obtain the
asymptotics
 of $V_0(s_1)$ as $s_1\to 0+$, we  use the standard theory
described, for example, in \cite{Bender}, Sec.3.4. To find the
leading terms of the asymptotic expansion for sufficiently small
$s_1>0$ we first take into account that \eqref{aid_cond2} implies
that there exists $\gamma=\rm const$ such that $f(\eta)\sim f_0
\eta^\gamma,$ $\eta\to \infty$, $f_0={\rm const}>0$, condition
\eqref{as_cond} implies that $\gamma\ge 1$. Thus, two linearly
independent solutions to \eqref{Vs1} behave  as
\begin{eqnarray}\label{Asg1}
Y_1&\sim & C\, s_1^\gamma,\quad Y_2\sim C\, \exp
\left(-\frac{\epsilon}{(\gamma-1)s_1^{\gamma-1}}\right),  \quad
\gamma>1,
\\
Y_1 &\sim& C\, s_1,\quad Y_2\sim C\, s_1^\epsilon,\quad \gamma=1,
\label{Ase1}
\end{eqnarray}
$s_1\to 0+$. Thus, if $\gamma>1$, we see from \eqref{Asg1} that
$(1-s)V(s,q_0(s))\sim \frac{V_0(s_1)}{s_1}= o(s_1)$, $s_1\to 0+$ or
$s\to -\infty$ and the behavior of $q_1(s)$, given by \eqref{q1}, is
defined only by term $\dfrac{f(1-s)}{1-s}$.

If $\gamma = 1$, then basically the last term in \eqref{q1} is
greater than  $\dfrac{f(1-s)}{1-s}$ as $s\to -\infty$, see
\eqref{Ase1}, and tends to plus or minus infinity depending on the
sign of $V$. The initial data can be chosen so that this term tends
to plus infinity and changes the behavior of the phase trajectory
\eqref {q1} in such a way that it remains in the lower half-plane,
and $ q $ tends to minus infinity, and $ s $ tends to minus
infinity.

As for the higher-order terms, $q_i$, $i=2,\dots$, they obey a
system of linear equations, similar to \eqref{al_k}, therefore for a
fixes $s$, the property to come back in the upper half-plane $q>0$
is defined only by $q_1(s)$ for sufficiently small $\epsilon$.
However at any $\epsilon>0$ for $s,q\to-\infty$ the last term in
\eqref{aid_qs} is subjected to the previous one, therefore, equation
\eqref{aid_qs} is equivalent at the point of singularity formation
to \eqref{as_phase}, and the phase trajectory for any initial data
turns out in the upper half-plane $q>0$, where $q,s$ cannot go to
infinity. Theorem \ref{aid_T} is proved. $\Box$

\begin{remark} We notice that for the non-perturbed case
$\epsilon=0$ the derivatives of solution can go to infinity at
$V_\infty=\lim\limits_{s=-\infty} V(s),$   this value is defined by
initial data and it can be any constant. For the case $\epsilon>0$,
the blow-up necessarily happens for $V_\infty=0$.
\end{remark}

\begin{remark} Note that for the case of affine solutions, the threshold friction $ \nu (n) = n $
guarantees the global smoothness of the solution, but for arbitrary data it is insufficient.
\end{remark}


\section{Discussion}\label{S4}

1. System \eqref{2} has the form $U_t + A(U,U_x) U_x =F(U) $,
$U=(V,E)$. The matrix $A$ is a Jordan block, it has  multiple
eigenvalue $V$, but only one eigenvector. The system does not belong
to symmetric hyperbolic one, and nonlinear resonance can occure in
the solution \cite{Isaacson_Temple}. The simplest system of this
form is the so-called pressureless
gas dynamics.   
It is commonly known that there the component of density develops
the delta-singularity. In our case it happens at the points where
$E_x$ tends to $-\infty$.

2. There exist different approaches to the well-posedness of weak
solutions to the pressureless Euler-Poisson equations
\cite{Isaacson_Temple}, \cite{LeFloch}, \cite{Tadmor_Wei}.

\section*{Acknowledgements}
The work is partially supported by the Moscow Center for Fundamental
and Applied Mathematics.
 The author thanks Evgenyi Chizhonkov, Alexander Kurganov and Maria Delova for a stimulating
 discussion.


\begin{thebibliography}{99}

\bibitem{Bender}   C.M. Bender, S.A. Orszag, Advanced mathematical methods for scientists and engineers,
 International series in pure and applied mathematics, McGraw-Hill,
 1978.

\bibitem{ELT} S. Engelberg,
H. Liu,  E. Tadmor. Critical thresholds in Euler-Poisson equations,
Indiana University Mathematics Journal, 50 (2001) 109-157. 



\bibitem{ChenWang} D. Wang, G.-Q. Chen,   Formation of singularities in compressible Euler-Poisson
fluids with heat diffusion and damping relaxation, Journal of
Differential Equations, 144 (1998) 44-65.

\bibitem{Ch_book} E.V. Chizhonkov, Mathematical Aspects of Modelling Oscillations and Wake Waves in Plasma, CRC Press, Boca Raton, 2019.

\bibitem{brodin}G. Brodin, L. Stenflo, Nonlinear dynamics of a cold collisional electron plasma, Physics of Plasmas 24(2017),
124505.

\bibitem{RCh_ZAMP} O.S. Rozanova, E.V. Chizhonkov,  On the conditions for the breaking of oscillations in a cold
plasma, Z. Angew. Math. Phys. 72(2021), 13,
\url{https://doi.org/10.1007/s00033-020-01440-3}

\bibitem{DRCh_AIP} O. Rozanova, E. Chizhonkov,
M. Delova, Exact thresholds in the dynamics of cold plasma with
electron-ion collisions, AIP Conference Proceedings, 2302 (2020),
060012, \url{https://doi.org/10.1063/5.0033619}

\bibitem{Isaacson_Temple} E. Isaacson, B. Temple,  Nonlinear resonance
in systems of conservation laws, Siam J. Appl. Math., 52 (1992)
1260-1278.


\bibitem{LeFloch}
P. LeFloch, S. Xiang,  Existence and uniqueness results for the
pressureless Euler-Poisson system in one spatial variable.
Portugaliae Mathematica, 72 (2015) 229-246.

\bibitem{Tadmor_Wei}
E. Tadmor, D. Wei, A variational representation of weak solutions
for the pressureless Euler-Poisson equations,     arXiv:1102.5579.



\end{thebibliography}
\end{document}